\documentclass[12pt]{amsart}

\usepackage{amssymb}
\usepackage{latexsym}
\usepackage{verbatim,euscript}
\usepackage{fullpage,color}

\input {cyracc.def}
\numberwithin{equation}{section}

\usepackage[colorlinks=true,linkcolor=blue,urlcolor=my_color,citecolor=magenta]{hyperref}

\definecolor{my_color}{rgb}{0,0.5,0.5}
\definecolor{MIXT}{rgb}{0.4,0.3,0.6}

 \font\tencyr=wncyr10 
\font\tencyi=wncyi10 
\font\tencysc=wncysc10 
\def\rus{\tencyr\cyracc}
\def\rusi{\tencyi\cyracc}
\def\rusc{\tencysc\cyracc}

\newtheorem{thm}{Theorem}[section] 

\newtheorem{utv}[thm]{Claim}
\newtheorem{lm}[thm]{Lemma}
\newtheorem{cor}[thm]{Corollary}
\newtheorem{prop}[thm]{Proposition}

\theoremstyle{remark}
\newtheorem{rmk}[thm]{Remark}

\theoremstyle{definition}

\newtheorem*{ex-bn}{Example}
\newtheorem*{rema}{Remark}

\newenvironment{proof*}
{\noindent {\sl Proof.}\quad }{\hfill
$\square$}

\newcommand {\ah}{{\mathfrak a}}
\newcommand {\be}{{\mathfrak b}}

\newcommand {\g}{{\mathfrak g}}
\newcommand {\h}{{\mathfrak h}}
\newcommand {\ka}{{\mathfrak k}}

\newcommand {\m}{{\mathfrak m}}

\newcommand {\q}{{\mathfrak q}}
\newcommand {\rr}{{\mathfrak r}}

\newcommand {\te}{{\mathfrak t}}
\newcommand {\ut}{{\mathfrak u}}

\newcommand {\slno}{\mathfrak{sl}_{l+1}}

\newcommand {\spn}{\mathfrak{sp}_{2l}}

\newcommand {\sono}{\mathfrak{so}_{2l+1}}
\newcommand {\sone}{\mathfrak{so}_{2l}}

\newcommand {\eus}{\EuScript}

\newcommand {\ap}{\alpha}

\newcommand {\N}{{\mathcal N}}
\newcommand {\co}{{\mathcal O}}

\newcommand {\BZ}{{\mathbb Z}}
\newcommand {\BN}{{\mathbb N}}

\newcommand {\md}{/\!\!/}

\newcommand {\Ad}{{\mathrm{Ad}}}

\newcommand {\codim}{{\mathrm{codim\,}}}

\newcommand {\gr}{{\mathrm{gr}}}

\newcommand {\ind}{{\mathsf{ind}}}
\newcommand {\Lie}{{\mathrm{Lie\,}}}

\newcommand {\Ima}{{\mathrm{Im\,}}}
\newcommand {\Mor}{\mathsf{Mor}}

\newcommand {\spe}{{\mathsf{Spec\,}}}

\newcommand {\trdeg}{{\mathrm{trdeg\,}}}

\newcommand {\GR}[2]{{\textrm{{\bf #1}}}_{#2}}

\newcommand {\ov}{\overline}

\newcommand {\beq}{\begin{equation}}
\newcommand {\eeq}{\end{equation}}

\renewcommand{\le}{\leqslant}
\renewcommand{\ge}{\geqslant}

\newcommand {\bbk}{\mathbb F}

\begin{document}
\setlength{\parskip}{2pt plus 4pt minus 0pt}
\hfill {\scriptsize July 4, 2011}
\vskip1.5ex

\title
{A remarkable contraction of semisimple Lie algebras}
\author[D.\,Panyushev]{Dmitri I.~Panyushev}
\address[D.P.]{Institute for Information Transmission Problems of the R.A.S., 
 B. Karetnyi per. 19, Moscow 127994, Russia}
\email{panyushev@iitp.ru}
\urladdr{\url{http://www.mccme.ru/~panyush}}
\author[O.\,Yakimova]{Oksana S.~Yakimova}
\address[O.Y.]{Universit\"at Erlangen-N\"urnberg,
Mathematisches Institut, Bismarckstra\ss e 1$\frac{1}{2}$, 
Erlangen D-91054, Deutschland}
\email{yakimova@mccme.ru}
\urladdr{\url{http://www.mi.uni-erlangen.de/~yakimova/}}
\maketitle

\section*{Introduction}
\noindent
The ground field $\bbk$ is algebraically closed and $\mathsf{char}\, \bbk=0$.
Let $G$ be a connected semisimple algebraic group of rank $l$ with  Lie algebra $\g$. 
Recently, E.\,Feigin introduced a very interesting contraction of $\g$  \cite{feigin1}. 
His motivation came from some problems in Representation Theory \cite{ffp}, and making use of this contraction he also studied certain degenerations of flag varieties  \cite{feigin2}.
Our goal is to elaborate on invariant-theoretic properties of 
these contractions of semisimple Lie algebras.

Fix a triangular decomposition
$  \g=\ut\oplus\te\oplus\ut^-$, where $\te$ is a Cartan subalgebra.
Then $\be=\ut\oplus\te$ is the fixed Borel subalgebra of $\g$. 
The corresponding subgroups of $G$ are $B,U$, and $T$.
Using the vector space isomorphism $\g/\be\simeq\ut^-$, we regard $\ut^-$ as a $B$-module. 
If $b\in\be$ and $\eta \in \ut^-$, then $(b,\eta) \mapsto b\circ \eta$ stands for the 
corresponding representation of $\be$. That is, if $p_-: \g\to \ut^-$ is the 
projection with kernel $\be$,  then $b\circ \eta=  p_-([b,\eta])$. 

Following \cite[Sect.\,2]{feigin1},  consider the semi-direct product
$\q=\be\ltimes (\g/\be)^a=\be\ltimes (\ut^-)^a$, where
the superscript `a' means that the $\be$-module $\ut^-$ is regarded as an abelian ideal in 
$\q$. 
We may (and will) identify the vector spaces $\g$ and $\q$ using the decomposition
$\g=\be\oplus\ut^-$.
If $(b,\eta), (b',\eta') \in \q$, then the Lie bracket in $\q$ is given by
\beq  \label{eq:skobka}
   [(b,\eta), (b',\eta')]=([b,b'], b\circ\eta'- b'\circ\eta) .
\eeq
The corresponding connected algebraic group is
$Q=B\ltimes N$, where $N=\exp((\ut^-)^a)$ is an abelian normal unipotent subgroup
of $Q$. 
The exponential map $\exp: (\ut^-)^a \to N$ is an isomorphism of varieties, and 
elements of $Q$ are written as product $s{\cdot}\exp(\eta)$, where $s\in B$ and $\eta\in\ut^-$.
If  $(s,\eta)\mapsto s{\centerdot}\eta$ is the representation of $B$ in $\ut^-$, then
the adjoint representation of $Q$ is given by
\beq   \label{eq:adj-action}
   \Ad_Q(s{\cdot}\exp(\eta))(b,\eta')=(\Ad(s)b, s{\centerdot} (\eta' - b\circ\eta)) .
\eeq
In this note, we explicitly construct certain polynomials that generate the algebras of 
invariants $\bbk[\q]^Q$ and $\bbk[\q^*]^Q$, and thereby prove that these two algebras are 
free. Furthermore, we also show that these polynomials generate the
corresponding fields of invariants, $\bbk(\q)^Q$ and $\bbk(\q^*)^Q$, and that
$\bbk[\q]$ is a free $\bbk[\q]^Q$-module and $\bbk[\q^*]$ is a free $\bbk[\q^*]^Q$-module.
The last assertion  implies
that the enveloping algebra of $\q$, $\eus U(\q)$, is a free module over its centre.
The Lie algebra $\q$ is a an {\it In\"on\"u-Wigner contraction\/} of $\g$ (see \cite[Ch.\,7 \S\,2.5]{t41}), and we also discuss the corresponding relationship between the invariants of $G$ and $Q$.

Certain classes of non-reductive algebraic Lie algebras $\q$ such that 
$\bbk[\q^*]^{Q}$ is a polynomial ring have been studied before. 
They include the centralisers of nilpotent elements in 
$\mathfrak{sl}_{l+1}$ and $\mathfrak{sp}_{2l}$ 
\cite{ppy}, 
$\BZ_2$-contractions of $\g$ \cite{coadj07}, and the truncated seaweed (biparabolic) 
subalgebras of $\mathfrak{sl}_{l+1}$ and $\mathfrak{sp}_{2l}$ \cite{jos}.
Our result enlarges this interesting family of  Lie algebras. 

Let $\q^*_{reg}$ denote the set of regular elements of $\q^*$, i.e., $x\in\q^*_{reg}$
if and only if $\dim Q{\cdot}x$ is maximal. For many problems related to coadjoint 
representations, it is vital to have that $\codim (\q^*\setminus\q^*_{reg})\ge 2$
 \cite{coadj07,ppy}.  However, we prove that if $\g$ is simple and not of 
type $\GR{A}{l}$, then $\q^*\setminus\q^*_{reg}$ contains a divisor.

\noindent
{\sl Notation.} \nopagebreak

-- \ the centraliser in $\g$ of $x\in\g$ is denoted by $\g^x$.  

-- \ $\varkappa$ is the Killing form on $\g$.

-- \ $\g_{reg}$ is the set of regular elements of $\g$, i.e., $x\in \g_{reg}$ if and only if 
$\dim\g^x=l$. 

-- \ If $X$ is an irreducible variety, then $\bbk[X]$ is the algebra of regular functions 
and $\bbk(X)$ is the field of rational functions on $X$. If $X$ is acted upon by an algebraic 
group $A$, then $\bbk[X]^A$ and $\bbk(X)^A$ denote the subsets of respective 
$A$-invariant functions.

-- \ If $\bbk[X]^A$ is finitely generated, then $X\md A:=\spe(\bbk[X]^A)$ and
$\pi: X\to X\md A$ is determined by the inclusion $\bbk[X]^A \hookrightarrow \bbk[X]$.
 If  $\bbk[X]^A$ is graded polynomial, then the elements of any set of algebraically independent homogeneous generators  will be referred to as {\it basic invariants\/}.

-- \ $\eus S^i(V)$ is the $i$-th symmetric power of the vector space $V$ and $\eus S(V)=\oplus_{i\ge 0}\eus S^i(V)$.

{\small {\bf Acknowledgments.} 
During the preparation of this paper, the second author benefited from an inspiring environment of the trimester  program ``On the Interaction of Representation Theory 
with Geometry and Combinatorics '' at HIM (Bonn).
She is grateful to P.\,Littelmann for the invitation.}

\section{On adjoint and coadjoint invariants of In\"on\"u-Wigner contractions} 
\label{sect0}

\noindent
The algebra $\q=\be\ltimes(\ut^-)^a$ is an In\"on\"u-Wigner contraction
of $\g$.
For this reason, we recall the relevant setting 
and then describe a general procedure for constructing adjoint and coadjoint invariants
of In\"on\"u-Wigner contractions.
The $\BZ_2$-{\it contractions\/} of $\g$ (considered in \cite{rims07,coadj07}) are special 
cases of  In\"on\"u-Wigner contractions, and for them such a procedure is exposed in 
\cite[Prop.\,3.1]{coadj07}. 
However, a more general situation considered here requires another proof.

For a while, we assume that $G$ is any connected algebraic group. Let $H$ be an arbitrary 
connected subgroup of $G$ and let $\m$ be a complementary subspace to $\h=\Lie H$. 
Using the vector space isomorphism $\g/\h\simeq\m$, we  regard $\m$ as $H$-module.
Consider the invertible linear map $\mathsf c_t: \g\to \g$, $t\in \bbk\setminus\{0\}$, 
such that $c_t(h+m)=h+tm$ \ ($h\in\h$, $m\in\m$) and
define the Lie algebra multiplication $[\ ,\ ]_{(t)}$
on the vector space $\g$ by the rule
\[
     [x,y]_{(t)}:= \mathsf c_t^{-1}\bigl( [\mathsf c_t(x), \mathsf c_t(y)]\bigr), \quad x,y\in\g \ .
\]
Write $\g_{(t)}$ for the corresponding Lie algebra.
The operator $(\mathsf c_t)^{-1}=\mathsf c_{t^{-1}}: \g\to \g_{(t)}$ yields an
isomorphism between the Lie algebras $\g=\g_{(1)}$ and $\g_{(t)}$,  hence
all algebras $\g_{(t)}$ are isomorphic.
It is easily seen that $\lim_{t\to 0}\g_{(t)}\simeq \h\ltimes (\g/\h)^a=\h\ltimes\m^a$. 

The resulting Lie algebra $\ka:=\h\ltimes\m^a$ is called an 
{\it In\"on\"u-Wigner contraction\/} of $\g$, cf. Example~7 in \cite[Chapter~7, \S\,2]{t41}.
The corresponding connected algebraic group is $K=H\ltimes \exp(\m^a)$. We identify the vector spaces $\g$ and  $\ka$ using the decomposition $\g=\h\oplus\m$.

\begin{rema} For $\g$ semisimple, the contraction $\g\leadsto \be\ltimes\ut^-$ is 
presented in a more lengthy way, using structure constants, in \cite[Remark~2.3]{feigin1}.
\end{rema}

\subsection{}   \label{subs:IW-inv}
To construct invariants of the coadjoint representation of $\ka$, we proceed as follows.
Let $f\in \eus S(\g)=\bbk[\g^*]$ be a homogeneous polynomial of degree $n$.
Using the  decomposition $\g=\h\oplus\m$, we consider the bi-homogeneous components of
$f$:
\[
     f=\sum_{a\le i\le b} f^{(n-i,i)} ,
\]
where $f^{(n-i,i)}\in \eus S^{n-i}(\h)\otimes \eus S^i(\m)\subset \eus S^n(\g)$, and
both $f^{(n-a,a)}$ and $f^{(n-b,b)}$ are assumed to be nonzero. In particular, 
$f^{(n-b,b)}$
is the bi-homogeneous component having the maximal degree relative to $\m$.
Since $\g_{(t)}$ and $\ka$ are just the same vector spaces, we also can regard each
$f^{(n-i,i)}$  as an element of $\eus S^n(\g_{(t)})$ or $\eus S^n(\ka)$.

\begin{thm}   \label{thm:coadj-inonu}
If  $f\in \eus S^n(\g)^G=\bbk[\g^*]^G_n$, then $f^{(n-b,b)}\in \eus S^n(\ka)^K=\bbk[\ka^*]^K_n$.
\end{thm}
\begin{proof}
The isomorphism of Lie algebras $\mathsf c_{t^{-1}}: \g\to \g_{(t)}$ implies that
$\sum_{a\le i\le b} t^{-i}  f^{(n-i,i)}\in  \eus S(\g_{(t)})^{G_{(t)}}$ for all $t\ne 0$. 
It is harmless to replace the last expression with the $G_{(t)}$-invariant
$f_{(t)}:=\sum_{a\le i\le b} t^{n-i}  f^{(n-i,i)}$. 
Consider the line $\langle f_{(t)}\rangle$ in the projective 
space $\mathbb P(\eus S^n(\g))$. It is a standard fact on $\bbk^\times$-actions in projective spaces that 
$\lim_{t\to 0} \langle f_{(t)}\rangle=\langle f^{(n-b,b)}\rangle$.
Hence $f^{(n-b,b)}$ is $K$-invariant.
\end{proof}

\noindent
Let us say that $f^\bullet:=f^{(n-b,b)}$ is the {\it highest component\/} of $f\in\bbk[\g^*]^G_n$
(with respect to the contraction $\g\leadsto \q$). 
Set $\eus L^\bullet(\bbk[\g^*]^G)=\{ f^\bullet \mid f\in\bbk[\g^*]^G \ \text{ is homogeneous}\}$.
Clearly, it 
is a graded subalgebra of $\bbk[\ka^*]^K$.

Invariants of the adjoint representation of $\ka$ can be constructed in a similar way.
Set $\m^*:=\h^\perp$, the annihilator of $\h$ in $\g^*$. Likewise, $\h^*=\m^\perp$.
Then $\g^*=\m^*\oplus\h^*$, and the adjoint operator
$\mathsf c_t^*: \g^*\to \g^*$ is given by
$\mathsf c_t^*(m^*+h^*)=t^{-1}m^*+h^*$ ($m^*\in\m^*$, $h^*\in\h^*$).
Having identified  $\q^*$ and $\ka^*$,
we can play the same game with homogeneous
elements of $\eus S(\g^*)=\bbk[\g]$.  If $\tilde f\in \eus S^n(\g^*)$, then $\tilde f^{(i,n-i)}$
denotes its bi-homogeneous component that belong to 
$\eus S^{i}(\m^*)\otimes \eus S^{n-i}(\h^*)$.
The resulting assertion is the following:

\begin{thm}  \label{thm:adj-inonu}
For $\tilde f\in \eus S^n(\g^*)^{G}$, let $\tilde f^{(a,n-a)}$ be the bi-homogeneous component
with minimal $a$, i.e., having the maximal degree relative to $\h^*=\m^\perp$.  
Then $\tilde f^{(a,n-a)}\in \eus S^n(\ka^*)^K$.
\end{thm}

Likewise, we write $\tilde f^\bullet:=\tilde f^{(a,n-a)}$ and
consider the algebra of highest components,  
$\eus L^\bullet(\bbk[\g]^G)$, which is a graded subalgebra of $\bbk[\ka]^K$.

\begin{lm}   \label{lm:poincare}
The graded algebras $\bbk[\g^*]^G$ and $\eus L^\bullet(\bbk[\g^*]^G)$ have the same Poincar\'e series, i.e., 
$\dim\bbk[\g^*]^G_n=\dim\eus L^\bullet(\bbk[\g^*]^G_n)$ for all $n\in\BN$; and likewise for\/ $\bbk[\g]^G$ and $\eus L^\bullet(\bbk[\g]^G)$.
\end{lm}
\begin{proof}
It easily seen that, for each $n$, there is a basis $\{p_i\}$ for $\bbk[\g^*]^G_n$ such that
$\{p_i^\bullet\}$ is a basis for $\eus L^\bullet(\bbk[\g^*]^G_n)$.
\end{proof}

It is not always the case that $\eus L^\bullet(\bbk[\g^*]^G)=\bbk[\ka^*]^K$ or
$\eus L^\bullet(\bbk[\g]^G)=\bbk[\ka]^K$.
For instance, we will see below that, for $\g$ semisimple and  $\q=\be\ltimes(\ut^-)^a$, such 
an equality holds only for the invariants of the coadjoint representation.
By the very construction, the algebras $\eus L^\bullet(\bbk[\g^*]^G)$ and $\eus L^\bullet(\bbk[\g]^G)$
are bi-graded. Moreover, it follows from \cite[Theorem\,2.7]{coadj07} that 
the algebras $\bbk[\ka^*]^{K}$ and $\bbk[\ka]^{K}$ are always bi-graded.

\subsection{}
If $\g$ is semisimple, then we may identify $\g$ and $\g^*$ (and hence
$\eus S(\g)$ and $\eus S(\g^*)$) using the Killing form
$\varkappa$. If $\h$ is also reductive, then  $\varkappa$ is non-degenerated on $\h$ and 
one can take $\m$ to be the orthocomplement of $\h$ with respect to $\varkappa$. 
Then $\h^\perp\simeq \m$ and the decompositions
of $\g$ and $\g^*$ considered in the general setting of In\"on\"u-Wigner contractions
coincide. Moreover, we can also identify the vector spaces $\ka$ and $\ka^*$. 
However, to obtain invariants of the adjoint and coadjoint
representations of $\q$, one has to take the bi-homogeneous components of maximal
degree with respect to {\it different\/} summands in the sum $\g=\h\oplus\m$.
In this situation, Theorems~\ref{thm:coadj-inonu} and \ref{thm:adj-inonu} admit
the following simultaneous formulation:

{\it Suppose that $f\in \bbk[\g]^G_n\simeq \eus S(\g)^G_n$  and 
$f=\sum_{a\le i\le b} f^{(n-i,i)}$ is the bi-homogeneous decomposition relative to the
sum $\g=\h\oplus\m$. (That is, $\deg_\h  f^{(n-i,i)}=n-i$, etc.)
Then, upon identifications of vector spaces $\g$,$\ka$, and $\ka^*$, we have
$f^{(n-a,a)}\in \bbk[\ka]^K$ and $f^{(n-b,b)}\in \bbk[\ka^*]^K$.}

Such a phenomenon was already observed in the case of $\BZ_2$-contractions of
semisimple Lie algebras, i.e., if $\h$ is the fixed-point subalgebra of an involution,
see \cite[Prop.\,3.1]{coadj07}.


\section{Invariants of the adjoint representation of $Q$} 
\label{sect1}

\noindent
In this section, we describe the algebra of invariants of the adjoint representation of
$Q$.

\noindent
To prove that a certain set of invariants generates  the whole
algebra of invariants, we use the following  lemma of  Igusa.

\begin{lm}[Igusa]   \label{lm:igusa}
Let $A$ be an algebraic group acting regularly on an irreducible affine variety 
$X$. Suppose that $S$ is an integrally closed finitely generated
subalgebra of\/ $\bbk[X]^A$ and the morphism $\pi: X\to \spe S=:Y$ has the properties:

{\sf (i)} \ the fibres of $\pi$ over a dense open subset of\/ $Y$ 
contain a dense $A$-orbit;

{\sf (ii)}  \ $\Ima\pi$ contains an open subset $\Omega$ of\/ $Y$ such that
$\codim (Y\setminus \Omega)\ge 2$.
\\
Then $S=\bbk[X]^A$. In particular, the algebra of $A$-invariants is finitely generated.
\end{lm}

See e.g. \cite[Lemma\,6.1]{rims07} for the proof.
\begin{rmk}  \label{rmk:zamena}
The proof given in \cite{rims07} shows that the above condition (i) can be replaced
with the condition that $S\subset \bbk[X]^A$ generates the field $\bbk(X)^A$.
(In fact, it is not hard to prove that (i) holds {\sl if and only if\/}  $S$ 
separates $A$-orbits in a dense open subset of $X$ {\sl if and only if\/} 
$S$ generates $\bbk(X)^A$.)
\end{rmk}

\begin{lm}  \label{lm:clear1}
If $t\in \te$ is regular and $u\in\ut$ is arbitrary, then
{\sf (i)} \  $t+u$ and $t$
belong to the same $\Ad\,U$-orbit;
{\sf (ii)} \  $(t+u)\circ \ut^-=\ut^-$.
\end{lm}
\begin{proof}
(i)  Clearly, $(\Ad\,U)t\subset t+\ut$ for all $t\in \te$. If $t$ is regular, then
$\dim (\Ad\,U)t=\dim \ut$. It is also known that the orbits of a unipotent group acting on an 
affine variety are closed. Hence $(\Ad\,U)t=t+\ut$.

(ii) This is obvious if $u=0$. In general, this follows from (i).
\end{proof}

\begin{thm}   \label{thm:adj}
We have 
$\bbk[\q]^Q\simeq \bbk[\te]$, and the quotient morphism $\pi_Q: \q\to \te$ is given by
$(t+u, \eta) \mapsto t$.
\end{thm}
\begin{proof}
Clearly, $\bbk[\q]^Q=(\bbk[\q]^{N})^B$. We  prove that \ 
1) $\bbk[\q]^{N}\simeq \bbk[\be]$ and \ 2)  $\bbk[\be]^B\simeq \bbk[\te]$.

1) \ Consider the projection $\pi_N: \q\to \q/(\ut^-)^a \simeq \be$. Clearly, $N$ acts trivially
on $\q/(\ut^-)^a$ and $\pi_N$ is a surjective $N$-equivariant morphism. Hence $\bbk[\be]\subset \bbk[\q]^N$. 
By Lemma~\ref{lm:igusa},
the equality $\bbk[\be]= \bbk[\q]^N$ will follow from the fact that
generic fibres of $\pi_N$ are $N$-orbits.

If $t\in\te$ is regular and $u\in\ut$ is arbitrary, then $b=t+u$ is a regular semisimple
element of $\g$. By \eqref{eq:adj-action} with $s=1$, we have
\[
   \Ad_Q(N)(b,\eta)=  (b,\eta+ b\circ \ut^-) .
\]
It then follows from Lemma~\ref{lm:clear1} that $\Ad_Q(N)(b,\eta)=(b,\ut^-)$.
On the other hand, $\pi_N^{-1}(b)=(b,\ut^-)$, i.e., $\pi_N^{-1}$ is a single $N$-orbit
whenever $b$ is regular semisimple.

2) \ Consider the projection $\pi_B: \be\to \be/\ut\simeq \te$. Clearly, $B$ acts trivially
on $\be/\ut$ and $\pi_B$ is a surjective $B$-equivariant morphism. Hence
$\bbk[\te]\subset \bbk[\be]^B$.
By Lemma~\ref{lm:igusa},
the equality $\bbk[\te]= \bbk[\be]^B$ will follow from the fact that
generic fibres of $\pi_B$ are $B$-orbits. 
Again, 
it follows from Lemma~\ref{lm:clear1} that if $t\in\te$ is regular, then $(\Ad\,B)t=t+\ut=
\pi_B^{-1}(t)$.
\end{proof}

\begin{rmk}     \label{rmk-adj}
Theorem~\ref{thm:adj} can be proved in a less informative way.
Notice that $[\q,\q]=\ut\ltimes(\ut^-)^{a}$ and therefore 
$\bbk[\mathfrak t]\subset\bbk[\q]^Q$. 
Let $x\in\te$ be regular semisimple. 
Then $\q^x\simeq \g^x=\te$, since $\g$ and $\q$ are isomorphic as $T$-modules. 
The fibres of the morphism $\pi_Q: \q\to \te$, defined in Theorem~\ref{thm:adj},
are linear spaces of 
dimension $\dim\q-\dim\te=\dim (\Ad\,Q)x$.
Hence a generic fibre contains a
dense $Q$-orbit and Lemma~\ref{lm:igusa} applies. 
We  also see that the algebra $\bbk[\te]$ separates generic $Q$-orbits in $\q$
and therefore $\bbk(\q)^Q=\bbk(\te)$. 
\end{rmk}

Comparing with the adjoint representation of $\g$, we see that, for $\q$, the algebra of 
invariants remains polynomial, but the degrees of basic invariants drastically decrease!
All the basic invariants in $\bbk[\q]^Q$ are of degree 1. This clearly means that here
$\eus L^\bullet(\bbk[\g]^G)\subsetneqq \bbk[\q]^Q$.

\section{Invariants of the coadjoint representation of $Q$} 
\label{sect2}

\noindent
In this section, we describe the algebra of invariants of the coadjoint representation 
of $Q$. 
The coadjoint representation is much more interesting since 
$\bbk[\q^*]=\eus S(\q)$ is a Poisson algebra, $\eus S(\q)^Q$ is the centre of this Poisson algebra, and 
$\eus S(\q)$ is related to the enveloping algebra of 
$\q$ via the Poincar\'e-Birkhoff-Witt theorem.

Since $\q$ is isomorphic to $\be\oplus \g/\be \simeq \be\oplus \ut^-$ as vector space, 
the dual vector space $\q^*$ is isomorphic to $(\g/\be)^*\oplus \be^*$. Using $\varkappa$, we 
identify $\be^*$ with $\be^-{:}=\te\oplus\ut^-$ and $(\g/\be)^*$ with $\ut$.
To stress that $\q^*$ is regarded as a $Q$-module and $\be^-$ appears to be a 
$Q$-stable subspace, we write
$\q^*=\ut{\,\propto\,}\be^-$.
If $(b,\eta)\in \q$ and $(u, \xi)\in \q^*$, i.e., $u\in\ut$ and $\xi\in\be^-$, 
then the coadjoint representation of $\q$ is given by the formula:
\beq   \label{eq:coadj}
    (b,\eta) {\star} (u, \xi)= ([b,u], \phi(u,\eta)+b\star\xi) .
\eeq
Here $(b,\xi) \mapsto b\star\xi$ is the coadjoint representation of $\be$, and
\[
\phi: \ut\times\ut^- \simeq \ut\times \ut^* \stackrel{\psi}{\to} \be^*\simeq \be^-,
\] 
where $\psi$ is  the {\it moment map\/} associated with the $\be$-module $\ut$.
Upon our identifications, the mapping $\phi$ is directly defined by
$\varkappa( b, \phi(u,\eta)):=\varkappa( [b,u], \eta)=- \varkappa( u, b\circ \eta)$.

Recall some well-known properties of the $B$-module $\ut$:

\textbullet \ \ 
If $\tilde e\in\ut $ is regular nilpotent, then $\g^{\tilde e}\subset \ut$ 
and hence $(\Ad\,B)\tilde e$ is dense in $\ut$. 

\textbullet \ \ For any $e\in \ut$, each irreducible component of $(\Ad\,G) e\cap\ut$
has dimension $\frac{1}{2}\dim (\Ad\,G) e$ \cite[4.3.11]{spr84}. 

\vskip1.3ex\noindent
Let $\Mor_G(\g,\g)$ denote the $\bbk[\g]^G$-module of polynomial $G$-equivariant
morphisms $F:\g\to \g$. By work of Kostant \cite{ko63},
$\Mor_G(\g,\g)$ is a free graded $\bbk[\g]^G$-module of rank $l$.
It was noticed by Th.\,Vust \cite[Char.\,III, \S\,2]{vust} (see also \cite{rais})
that  a homogeneous basis of this module is obtained as follows.
Let $f_1,\dots,f_l$ be homogeneous algebraically independent generators of $\bbk[\g]^G$.
Each differential $\textsl{d}f_i$ determines a polynomial $G$-equivariant morphism
(covariant) from $\g$ to $\g^*$.
Identifying $\g$ with $\g^*$ via $\kappa$ yields a homogeneous covariant
(or, vector field) $F_i={\mathsf{grad}}f_i : \g\to \g$. Then 
$F_1,\dots,F_l$ form a homogeneous basis for 
$\Mor_G(\g,\g)$.
If $\deg f_i=d_i$, then $\deg F_i=d_i-1=:m_i$. It is customary to say that $\{m_1,\dots,m_l\}$
are the {\it exponents\/} of (the Weyl group of) $\g$. Recall that if $\g$ is simple and
$m_1\le \dots\le m_l$, then $m_1=1$, $m_2\ge 2$, and 
$m_i+m_{l-i+1}$ is the Coxeter number of $\g$.

The covariants $F_i$ have the following properties: 

(i)  $F_i(x)\in \g^x$ for all $i\in\{1,2,\dots,l\}$ and $x\in \g$;

(ii) The vectors $F_1(x),\dots,F_l(x)\in\g$ are linearly independent if and only if $x\in \g_{reg}$
\cite[Theorem\,9]{ko63}.

\noindent
It follows that $(F_1(x),\dots,F_l(x))$ is a basis for $\g^x$ if and only if $x\in \g_{reg}$.

\begin{lm}  \label{lm:covar-B}
If $x\in \be$, then $F_i(x)\in \be$. 
If $y\in \ut$, then $F_i(y)\in \ut$. 
\end{lm}
\begin{proof}
If $x\in\be\cap\g_{reg}$, then $\g^x\subset \be$. Hence $F_i(x)\in \g^x\subset\be$.
Since $\be\cap\g_{reg}$ is open and dense in $\be$, the assertion follows.
\\
If $y\in\ut\cap\g_{reg}$, i.e., $y$ is regular nilpotent, then $\g^y\subset \ut$ \cite{ko63}.
The rest is the same.
\end{proof}
Consequently, letting $P_i:=F_i\vert_\ut$, we obtain the covariants $P_1,\dots,
P_l \in \Mor_B(\ut,\ut)$. Actually, we consider the $P_i$'s as $B$-equivariant morphisms 
$P_i:\ut\to \ut\subset \be$.
Using these covariants, we define polynomials 
$\widehat P_i\in \bbk[\q^*]=\bbk[\ut{\,\propto\,}\be^-]$ by the formula
\beq  \label{eq:inv-coadj}
    \widehat P_i(u,\xi)= \varkappa( P_i(u),\xi ) , \quad i=1,\dots,l ,
\eeq
where $u\in\ut$ and $\xi\in\be^-$.
\begin{lm}  \label{lm:Q-inv}
We have $\widehat P_i \in \bbk[\q^*]^Q$.
\end{lm}
\begin{proof}
Since $Q=B\ltimes N$, it suffices to verify that $\widehat P_i$ is both $B$- and $N$-invariant.

1)  $\widehat P_i$ is $B$-invariant, since $P_i$ is $B$-equivariant.

2) 
For polynomials obtained from covariants 
$P_i$ as in \eqref{eq:inv-coadj}, 
the invariance with respect to the commutative unipotent group $N$ is equivalent to that
$[P_i(u),u]=0$. Indeed, for $\eta\in\ut^-$, the coadjoint action of $\exp(\eta)\in N$ is given by
$\exp(\eta){\star}(u,\xi)=(u,\xi+ \phi(u,\eta))$. Then
\begin{multline*}
 \widehat P_i(\exp(\eta){\cdot}(u,\xi))=\varkappa( P_i(u), \xi+ \phi(u,\eta))=
\varkappa( P_i(u), \xi)) +\varkappa(P_i(u), \phi(u,\eta))= \\
  \widehat P_i(u,\xi)+ \varkappa( [P_i(u),u], \eta) .
\end{multline*}
Hence $ \widehat P_i(\exp(\eta){\cdot}(u,\xi))=\widehat P_i(u,\xi)$  for all $\eta$ if and only if
$[P_i(u),u]=0$. In our case, this follows from the corresponding property of $F_i$. 
\end{proof}

\begin{rema}
We prove below that $\widehat P_i$ is the highest components of $f_i \in \bbk[\g^*]^G$.
In view of Theorem~\ref{thm:coadj-inonu}, this also implies that $\widehat P_i$ is
$Q$-invariant.
\end{rema}

\begin{thm}  \label{thm:coadj1}
The algebra $\bbk[\q^*]^Q$ is freely generated by $\widehat P_1,\dots, \widehat P_l$, 
and\/  $\bbk(\q^*)^Q$ is the fraction field of\/ $\bbk[\q^*]^Q$.
\end{thm}
\begin{proof}
Consider the morphism
\[
   \pi:  \q^*=\ut{\,\propto\,}\be^- \to \mathbb A^l ,
\]
given by $\pi(u,\xi)=(\widehat P_1(u,\xi),\dots,\widehat P_l(u,\xi))$. As in Section~\ref{sect1}, 
to prove that $\pi$ is the quotient by $Q$, we are going to apply
Lemma~\ref{lm:igusa} to $\pi$. 

If $e\in\ut$ is regular, then $P_1(e),\dots,P_l(e)$ are linearly independent and form
a basis for $\g^e=\ut^e$. Therefore, \eqref{eq:inv-coadj} implies that $\pi$ is onto, 
and condition~(ii) in Lemma~\ref{lm:igusa} is satisfied.

Let us prove that $\bbk(\q^*)^Q=\bbk(\widehat P_1,\dots, \widehat P_l)$.
Consider the morphism
\[
   \tilde\pi: \q^*\to (\q^*/\be^-)\times \mathbb A^l=\ut \times \mathbb A^l
\]
defined by 
$\tilde\pi(u,\xi)=(u,\widehat P_1(u,\xi),\dots,\widehat P_l(u,\xi))$.
If $e\in \ut\cap\g_{reg}$, then Eq.~\eqref{eq:inv-coadj} shows that
$\tilde\pi^{-1}(e, {a})$ is an affine subspace of $\q^*$ for any ${a}\in \mathbb A^l$, and
$\dim\tilde\pi^{-1}(e,{a})=\dim\be-l=\dim \ut$.
As in the proof of Theorem~\ref{thm:adj}, this implies that $\tilde\pi^{-1}(e,{a})$
is a sole $N$-orbit. Thus, the coordinate functions on $\ut$ and 
$\widehat P_1,\dots, \widehat P_l$ separate
generic $N$-orbits of maximal dimension. By the Rosenlicht theorem, this implies
that all these functions  generate the field of $N$-invariants on $\q^*$, i.e.,
$\bbk(\q^*)^N=\bbk(\ut)(\widehat P_1,\dots, \widehat P_l)$.
Since $B$ has an open orbit in $\ut$, we have $\bbk(\ut)^B=\bbk$.
Hence
\[
  \bbk(\q^*)^Q=(\bbk(\ut)(\widehat P_1,\dots, \widehat P_l))^B=\bbk(\widehat P_1,\dots, \widehat P_l) .
\]
In view of Remark~\ref{rmk:zamena}, this is sufficient for using Lemma~\ref{lm:igusa},
and we conclude that $\widehat P_1,\dots, \widehat P_l$ generate the algebra of 
$Q$-invariants on $\q^*$.
\end{proof}

\begin{rmk}  \label{rmk:new-inv}
Although we have proved that $\bbk(\q^*)^N=\bbk(\ut)(\widehat P_1,\dots, \widehat P_l)$,
it is not true that $\bbk[\q^*]^N=\bbk[\ut][\widehat P_1,\dots, \widehat P_l]$.
The reason is that the morphism $\tilde\pi$ defined in the previous proof does not 
satisfy condition (ii) of Lemma~\ref{lm:igusa}. That is, the closure of the complement
of $\Ima\tilde\pi$ contains a divisor.
One can prove that if $D=\ut\setminus (\Ad\,B)\tilde e=\ut\setminus (\ut\cap\g_{reg})$, then
this divisor is $D\times \mathbb A^l$.
Actually, we can explicitly point out a function in
$\bbk[\q^*]^N\setminus\bbk[\ut][\widehat P_1,\dots, \widehat P_l]$.
Let $v$ be a non-zero vector in the one-dimensional space $\be^U$. We can regard $v$ 
as a linear function on $\be^-$ and hence on $\q^*$. Making use of Eq.~\eqref{eq:skobka},
it is not hard to check that the subalgebra $(\ut^-)^a\subset \q$ commutes with $v$, i.e., 
$v$ is a required $N$-invariant in the symmetric algebra $\eus S(\q)$.
\end{rmk}

Recall that, for an algebraic group $A$ with Lie algebra $\ah$,
the {\it index} of  $\ah$, $\ind\,\ah$, equals 
$\trdeg \bbk(\ah^*)^A$. It is also true that the index cannot decrease under contractions,
hence $\ind\,\q\ge \ind\,\g=l$. The above description of the field of $Q$-invariants implies that
 
\begin{cor}
$\ind\,\q=l$.
\end{cor}
\begin{thm}   \label{thm:coadj2}
The polynomial ring $\bbk[\q^*]$ is a free $\bbk[\q^*]^Q$-module.
\end{thm}
\begin{proof}
Since it is already known that $\bbk[\q^*]^Q$ is a polynomial algebra (of Krull dimension $l$),
it suffices to prove that the quotient morphism 
$\pi: \q^*\to \q^*\md Q\simeq \mathbb A^l$ is equidimensional
\cite[Prop.\,17.29]{lift}. This, in turn, will follow from the fact that
the null-cone $\N=\pi^{-1}(\pi(0))$ is of dimension $\dim\q -l$.
To estimate the dimension of $\N$, consider the projection $p:\N\to \ut$ and partition
$\ut$ into finitely many {\it orbital varieties}, i.e., the irreducible components
of $(\Ad\,G)e_i\cap\ut$,  where $\{e_i\}$ runs over a finite set of representatives of all nilpotent 
$G$-orbits. Let $Z_i$ be an irreducible component of $(\Ad\,G)e_i\cap\ut$. Since
$\pi=(\widehat P_1,\dots, \widehat P_l)$, Eq.~\eqref{eq:inv-coadj}
shows that
\[
\dim p^{-1}(Z_i)=\dim Z_i+ \dim\be- \dim \text{span}\{P_1(e_i)\dots, P_l(e_i)\} .
\]
As $\dim Z_i=\frac{1}{2}\dim (\Ad\,G)e_i$, the condition that $\dim p^{-1}(Z_i)\le \dim\q-l$
can easily be transformed into 
\beq     \label{eq:brilliant}
    \dim\g^{e_i} + 2\dim \text{span}\{P_1(e_i)\dots, P_l(e_i)\} \ge 3l .
\eeq
Recall that $P_1,\dots,P_l$ are just the restrictions to $\ut$ of basic covariants
$F_1,\dots,F_l$, and $F_j=\mathsf{grad}f_j$. Consequently, 
$\dim \text{span}\{P_1(e_i)\dots, P_l(e_i)\}$ equals the rank of the differential at $e$
of the quotient morphism
$\pi_{\g,G}: \g \to \g\md G$. Therefore,  \eqref{eq:brilliant} is preisely the inequality
proved in \cite[Theorem\,10.6]{rims07}.
\end{proof}

\begin{cor}
The enveloping algebra $\eus U(\q)$ is a free module over its centre $\eus Z(\q)$.
\end{cor}
\begin{proof}
This is a standard consequence of the fact that $\bbk[\q^*]=\eus S(\q)$ is a free module over
$\eus S(\q)^Q$,  $\eus S(\q)^Q$ is the centre of the Poisson algebra $\eus S(\q)$,
and $\gr\,\eus Z(\q)=\eus S(\q)^Q$,
cf. \cite[Theorem\,21]{ko63}, \cite[Theorem\,3.3]{geof}.
\end{proof}

\begin{rmk}
By Theorem~\ref{thm:coadj2}, the irreducible components of all fibres of
$\pi:  \q^* \to \q^*\md Q\simeq \mathbb A^l$ are of dimension $\dim\q-l$. However, unlike the case of the (co)adjoint representation of $\g$,
the zero fibre of $\pi$ is highly reducible. For, if
$\dim\g^{e_i} + 2\dim \text{span}\{P_1(e_i)\dots, P_l(e_i)\} = 3l$, then every irreducible component of $(\Ad\,G)e_i\cap\ut$ gives rise to an irreducible component of 
$\pi^{-1}(\pi(0))$. A complete classification of nilpotent elements of $\g$ satisfying this equality 
is contained in \cite[\S\,10]{rims07}.
\end{rmk}

\begin{thm}   \label{thm:coadj-highest}
We have $\eus L^\bullet(\eus S(\g)^G)=\eus S(\q)^Q$.
The polynomials $\widehat P_1,\dots, \widehat P_l\in \bbk[\q^*]^Q=\eus S(\q)^Q$ are 
the highest components of  $f_1,\dots, f_l\in \eus S(\g)^G$ in the sense of 
Subsection~\ref{subs:IW-inv}. 
\end{thm} 
\begin{proof}
1) \ Since $\deg \widehat P_i=\deg f_i$ for all $i$, it follows from Lemma~\ref{lm:poincare}
and Theorem~\ref{thm:coadj1}
that $\eus L^\bullet(\eus S(\g)^G)$ and $\eus S(\q)^Q$ have the same Poincar\'e series. 
Hence these algebras coincide.

2) Recall that $\deg f_i=d_i=m_i+1$.
According to Theorem~\ref{thm:coadj-inonu}, we have to take the decomposition 
$\g=\be\oplus\ut^-$ 
and pick the bi-homogeneous component of $f_i$ of maximal degree with respect to $\ut^-$.

If the component $f_i^{(0,d_i)}\in\eus S^{d_i}(\ut^-)$ were non-trivial, then 
it would be a $Q$-invariant in $\eus S(\q)$ and in particular a $B$-invariant (Theorem~\ref{thm:coadj-inonu}).
Recall that if we work in $\q$, then $\ut^-\simeq \g/\be$ as $B$-module.
Since $\eus S(\g/\be)\simeq \bbk[\ut]$ and $\bbk[\ut]^B=\bbk$, we get a contradiction.
Hence $f_i^{(0,d_i)}=0$.

Then next possible component is $f_i^{(1,m_i)}\in \be\otimes \eus S^{m_i}(\ut^-)$. 
Using the identifications $\be^*\simeq \be^-$ and $\ut^*\simeq \ut^-$, 
we have $f_i^{(1,m_i)}\in \bbk[\be^-]_1\otimes \bbk[\ut]_{m_i}$. That is, 
if considered as a function on $\g=\be^-\oplus \ut$, it can be written as
$f_i^{(1,m_i)}(\xi,u)=\kappa( \bar P_i(u),\xi)$ for some morphism $\bar P_i: \ut \to \be$ of 
degree $m_i$. As we have already proved that $f_i^{(0,d_i)}=0$, $\bar P_i(u)$
is nothing but the value of $\mathsf{grad}f_i$ at $u$. Hence $\bar P_i=P_i$, and we are done.
\end{proof}

\section{Further properties of the coadjoint representation} 
\label{sect3}

\subsection{}   \label{subs:classical-covar}
For the classical Lie algebras, the basic covariants $F_i: \g\to \g$ (and hence $P_i$) have a 
simple description:
\begin{itemize}
\item if $x\in\slno$, then $F_i(x)=x^i$, $i=1,2,\dots,l$;
\item if $x\in\spn$ or $\sono$, then $F_i(x)=x^{2i-1}$, $i=1,2,\dots,l$;
\item if $x\in\sone$, then $F_i(x)=x^{2i-1}$, $i=1,2,\dots,l-1$. The covariant $F_l$ that is 
related to the pfaffian is described as follows. Let $x$ be a skew-symmetric
matrix  of order $2l$. For $i\ne j$, let $x_{[ij]}$ be the skew-symmetric sub-matrix 
of order $2l-2$
obtained by deleting $i$th and $j$th row and column. Set $a_{ij}=\mathsf{Pf}(x_{[ij]})$ if $i\ne j$, and
$a_{ii}=0$.  Then $F_l(x)=(a_{ij})_{i,j=1}^{2l}$. Clearly, $\deg F_l=l-1$, as required.
\end{itemize}

\noindent
Results of Sections~\ref{sect1} and \ref{sect2} explicitly yield the bi-degrees 
of basic invariants for $\q=\be\ltimes(\ut^-)^a$. 
For $\bbk[\q]^Q$, all the basic invariants have bi-degrees $(1,0)$. For $\bbk[\q^*]^Q$,  the basic invariants have bi-degrees $(m_i,1)$, i.e., 
$\widehat P_i\in \eus S^{m_i}(\ut^-)\otimes \be$. 
In particular, for the coadjoint representation, 
the total degrees of the basic $Q$-invariants remain the same as for $G$.

\subsection{}  
Hereafter we assume that $\g$ is simple and the  basic invariants $f_1,\dots,f_l\in
\bbk[\g]^G$ are numbered such that $d_i\le d_{i+1}$.  
Then $d_l=\mathsf h$ is the Coxeter number of $\g$.
We show that the corresponding
$Q$-invariant $\widehat P_l$ has a rather simple form. In fact, it appears to be a monomial.

Let $\Delta$ be the root system of $(\g,\te)$ and
$\Delta^+$ the subset of positive roots corresponding to $\ut$. 
Then $\Pi=\{\ap_1,\dots,\ap_l\}$ (resp. $\theta$) is the  set of simple roots 
(resp. the highest root) in $\Delta^+$.  Then $\theta=\sum_{i=1}^l a_i\ap_i$ and
$\sum_{i=1}^l a_i=\mathsf h-1$.
For any $\gamma\in\Delta$, $\g_\gamma$ denotes the corresponding root subspace,
and we fix a nonzero vector $e_\gamma\in\g_\gamma$. 

\begin{lm}    \label{lm:P_l}
Up to a scalar multiple, we have
$\widehat P_l= e_{-\ap_1}^{a_1}\dots e_{-\ap_l}^{a_l}e_\theta\in \eus S(\q)^Q$.
\end{lm}
\begin{proof} 
Recall that $\q=\be\oplus \ut^-$ as vector space, and  here $e_\theta\in \be$ and
$e_{-\ap_i}\in \ut^-$. By the very construction, $\widehat P:=
e_{-\ap_1}^{a_1}\dots e_{-\ap_l}^{a_l}e_\theta$ is a $T$-invariant in $\eus S(\q)$.
Then, using Eq.~\eqref{eq:skobka}, one readily verifies that $\widehat P$ is
both $U$-invariant and $N$-invariant.
Hence $\widehat P$ is a polynomial in $\widehat P_1,\dots,\widehat P_l$.
Since $\textsf{bi-deg\,} \widehat P=(\mathsf h-1,1)$ and  $m_i < m_l$ for $i<l$, 
the subspace of bi-degree
$(m_l,1)=(\mathsf h-1,1)$ in  $\eus S(\q)^Q$ is one-dimensional and spanned by
$\widehat P_l$. Hence the assertion.
\end{proof}

Since $\dim\q=\dim\g$, $\ind\,\q=\ind\,\g$, and the (total) degrees of the basic invariants of the coadjoint representations for $G$ and $Q$ coincide, we have the equality
\beq    \label{eq:arifm-tozhd}
   \sum_{i=1}^l \deg \widehat P_i=\frac{\dim\q+ \ind\,\q}{2} \ ,
\eeq
which is very useful in the study of the coadjoint representation, 
see e.g. \cite[Theorem\,1.2]{coadj07}. 
Unfortunately, $\q$ does not always possess another important ingredient, the so-called 
{\sl codim-2\/} property.
Recall that $x\in \q^*$ is said to be {\it regular\/} if $\dim Q{\cdot}x$ is maximal. The set of all 
regular elements is denoted by $\q^*_{reg}$. It is an open subset of $\q^*$,  and we say 
that $\q$ has the {\sl codim-2\/} property if $\codim(\q^*\setminus\q^*_{reg})\ge 2$.

\begin{thm}    \label{thm:no-codim2}
The algebra $\q$ does not have the {\sl codim-2\/} property if\/ $\g$ is not
of type $\GR{A}{l}$.
\end{thm}
\begin{proof}
Suppose that  $\q$ has the {\sl codim-2\/} property. Since \eqref{eq:arifm-tozhd} 
is satisfied, it follows from \cite[Theorem\,1.2]{coadj07} that 
the differentials $(\textsl d\widehat P_i)_x$, $i=1,\dots,l$, are linearly independent
if and only if $x\in \q^*_{reg}$. In particular, any divisor  $\tilde D\subset\q^*$  contains a point 
where the differentials of $\widehat P_1,\dots,\widehat P_l$ are linearly independent.

On the other hand, Lemma~\ref{lm:P_l} shows that if $a_i\ge 2$ for some $i$, then
$\textsl d\widehat P_l$ vanishes at the hyperplane  $\{ e_{-\ap_i}=0\}$, where $e_{-\ap_i}$
is regarded as a linear function on $\ut$ and hence on $\q^*$. Thus, $\q$ cannot have
the {\sl codim-2\/} property unless $a_i=1$ for all $i$, i.e.,  $\g$ is of type
$\GR{A}{l}$.
\end{proof}

To prove the converse of this theorem, we need some preparations.  For $\ap_i\in\Pi$,
let $\ut_i\subset \ut$ denote the kernel of the linear form $u\mapsto \kappa(e_{-\ap_i},u)$.
By \cite{ko63}, $\ut\setminus \ut\cap\g_{reg}=\cup_{i} \ut_i$.
Set 
\beq  \label{eq:Y}
 \mathcal Y= \mathcal Y(\q^*)=\{ x\in\q^* \mid  (\textsl d\widehat P_1)_x, \dots, (\textsl d\widehat P_l)_x 
 \ \text{ \ are linearly independent}\} .
 \eeq
 \begin{prop}   \label{prop:codim(Y)}
If \ $\g=\slno$, then $\codim(\q^*\setminus \mathcal Y)\ge 2$. 
\end{prop}
\begin{proof}
Let $a=(e,\xi)$ and $a'=(e',\xi')$
be typical elements of $\q^*$, where $e,e'\in\ut$ and $\xi,\xi'\in \be^-$.
According to formulae of Subsection~\ref{subs:classical-covar}, 
$\widehat P_i(e,\xi)=\kappa(e^i,\xi)$. Recall that $(\textsl{d}\widehat P_i)_a\in\q$
and $\langle (\textsl{d}\widehat P_i)_a, a'\rangle$ is the coefficient of $t$
in the expansion of $\widehat P_i(a+ta')$. Consequently, 
\[
    \langle (\textsl{d}\widehat P_i)_a, a'\rangle=\kappa(e^i,\xi')+\kappa(\sum_{k+m=i-1}
    e^ke' e^m, \xi) .
\]
The vector $(\textsl{d}\widehat P_i)_a$ has the $\be$- and $\ut^-$-components, and  this equality shows that:

\textbullet\quad the $\be$-component of $(\textsl{d}\widehat P_i)_a$ \ equals \ $e^i$;

\textbullet\quad the $\ut^-$-component of $(\textsl{d}\widehat P_i)_a$, say
$(\textsl{d}\widehat P_i)_a\{\ut^-\}$,  is determined by the equation
$\kappa((\textsl{d}\widehat P_i)_a\{\ut^-\}, e')=\kappa(\sum_{k+m=i-1} e^ke' e^m, \xi)$.

\noindent
Let $\co^{reg}$ and $\co^{sub}$ denote the regular and subregular nilpotent orbits in 
$\slno$, respectively.  Then $\ov{\co^{sub}\cap\ut} =\cup_j \ut_j$.
If $e\in\co^{reg}\cap\ut$, then
the $\be$-components of $(\textsl{d}\widehat P_i)_{(e,\xi)}$, $i=1,\dots,l$,
are linearly independent, regardless of $\xi$. Hence $(\co^{reg}\cap\ut)\times \be^- \subset \mathcal Y$.

If $e\in\co^{sub}\cap\ut$, then the $\be$-components of $(\textsl{d}\widehat P_i)_{(e,\xi)}$, 
$i=1,\dots,l-1$, are still linearly independent for any $\xi$, 
but $e^l=0$. However,  if $e$ is sufficiently general, then
the $\ut^-$-component of $(\textsl{d}\widehat P_l)_{(e,\xi)}$ appears to be nonzero 
for all $\xi$ that belong to a dense open subset of $\be^-$. More precisely, suppose that 
$e\in\ut_j$ and $\kappa(e,e_{-\ap_i}) \ne 0$   for $i\ne j$. Taking
$e'=e_{\ap_j}$, one readily computes that
$\sum_{k+m=l-1} e^ke' e^m=e^{j-1}e_{\ap_i}e^{l-j}$ is a nonzero multiple of $e_\theta$.
Hence, one can take any $\xi$ such that  $\kappa(\xi,e_{\theta})\ne 0$.

Thus, there is a dense open subset $\Omega\subset \cup_i \ut_i \times \be^-$
such that  $\Omega\subset \mathcal Y$, and the assertion follows.
\end{proof}

It turns out that Proposition~\ref{prop:codim(Y)} together with \eqref{eq:arifm-tozhd}
is sufficient to conclude that for $\g=\slno$,  $\q$ has the {\sl codim-2\/} property.
This follows from the following general assertion:

\begin{thm}   \label{thm:general-assert}
Let $R$ be a connected algebraic group with Lie algebra $\rr$.
Suppose that (a)  $\ind\,\rr=m$,  (b) $\bbk[\rr^*]^R=\bbk[p_1,\dots,p_m]$ is a graded 
polynomial algebra,
and  (c) $\sum_{i=1}^m \deg p_i=(\dim\rr+\ind\,\rr)/2$. Then the following conditions are 
equivalent:

1) \ $\codim(\rr^*\setminus\rr^*_{reg})\ge 2$;

2) \ $\codim(\rr^*\setminus \mathcal Y(\rr^*))\ge 2$, where $\mathcal Y(\rr^*)$
is defined as in \eqref{eq:Y} via the $p_i$'s.
\\ If these conditions are satisfied, then actually \ $\rr^*_{reg}=\mathcal Y(\rr^*)$.
\end{thm}
\begin{proof}
The implication $1) \Rightarrow 2)$ is already proved in \cite[Theorem\,1.2]{coadj07}.
\\
To prove the converse, one can slightly adjust the proof given in \cite{coadj07}, see also the 
proof of Theorem\,1.2 in \cite{ppy}.
Set $n=\dim\rr$. Let $T(\rr^*)$ denote the tangent bundle of $\rr^*$.
The main part of that proof consists in a construction of two homogeneous polynomial
sections of $\wedge^{n-m} T(\rr^*)$, denoted $\mathfrak V_1$ and $\mathfrak V_2$.
Write $(\mathfrak V_i)_x$ for the value of $\mathfrak V_i$ at $x\in\rr^*$.
These sections have the following properties:

{\bf --} \ There exist $F_1,F_2\in \bbk[\rr^*]$ such that $F_1\mathfrak V_1=F_2\mathfrak V_2$;

{\bf --} \ $(\mathfrak V_1)_x\ne 0$ if and only if $x\in\rr^*_{reg}$

{\bf --} \ $(\mathfrak V_2)_x\ne 0$ if and only if $x\in \mathcal Y(\rr^*)$;

{\bf --} \  $\deg \mathfrak V_1=(n-m)/2$ and $\deg\mathfrak V_2=\sum_i (\deg p_i-1)$.

\noindent
This only requires assumptions a) and b). If c) is also satisfied, then 
$\deg \mathfrak V_1=\deg \mathfrak V_2$.  Therefore either of conditions 1), 2)  implies
the other.
\end{proof}

Since $\q=\be\ltimes\ut^-$ does not have the {\sl codim-2\/} property if
$\g$ is not of type $\GR{A}{l}$,  we cannot  immediately conclude 
that in all cases  $x\in\q^*_{reg}$ if and only if 
$ (\textsl d\widehat P_1)_x, \dots, (\textsl d\widehat P_l)_x$  are linearly independent.
Nevertheless, the fact that 
$\widehat P_1,\dots,\widehat P_l$ are the highest components of the basic 
$G$-invariants $f_1,\dots,f_l$  allows to circumvent this difficulty.
It can be shown in general (see \cite{sim-z2}) 
that the coadjoint representation $(Q:\q^*)$ has the following property:

\begin{utv}
For  $x\in\q^*$ the following conditions are equivalent:
 \\
 \textbullet \quad The orbit $Q{\cdot}x$ is of maximal dimension, which is $\dim\q-l$ in this situation;
\\ 
 \textbullet \quad The differentials $(\textsl d\widehat P_i)_x$, $i=1,\dots,l$, are linearly independent.
\end{utv}
In case of the coadjoint representation of $\g$, this is a  result of Kostant
\cite[Theorem\,9]{ko63}.

\end{document}